\documentclass[12pt]{amsart}
\usepackage{amssymb,latexsym}
\usepackage{tikz}
\usetikzlibrary{positioning}
\usepackage{color}
\definecolor{Blue}{rgb}{0.3,0.3,0.9}
\definecolor{Red}{rgb}{.1,.1,.1}
\theoremstyle{plain}
\newtheorem{theorem}{Theorem}
\newtheorem{corollary}{Corollary}
\newtheorem*{main}{Main~Theorem}
\newtheorem*{mainc}{Main~Corollary}
\newtheorem{lemma}{Lemma}
\newtheorem{proposition}{Proposition}

\theoremstyle{definition}
\newtheorem{definition}{Definition}

\theoremstyle{definition}
\newtheorem{example}{Example}[section]

\theoremstyle{remark}

\newtheorem*{note}{Note}

\numberwithin{equation}{section}

\begin{document}
\title[Rough Set Algebra and CRDSA]{A Note on Rough Set Algebra and Core Regular Double Stone Algebras (CRDSA)}
\author[Daniel J. Clouse\\U.S. DoD - LACR\\USA]{Daniel J. Clouse\\U.S. DoD - Laboratory for Advanced Cybersecurity Research (LACR)\\USA}
\address{Daniel J. Clouse\\U.S. DoD - LACR\\USA} 
\email{djclouse@gmail.com}
\urladdr{}
\thanks{*Our thanks to referee Feng for suggesting we consider Quantum mechanics} 
\date{December 26, 2019}
\begin{abstract}
{\color{Red} Rough Set Theory (RST), first introduced by Pawlak in 1982, is an approach for dealing with information systems where knowledge is uncertain or incomplete.\cite{Pawlak} It is of fundamental importance in many subfields of artificial intelligence and cognitive science.\cite{RSTppf}} Given a universe $U$ with an equivalence relation $\theta$, the pair $\langle U,\theta\rangle$ is referred to as an information system and we denote its collection of rough sets $R_\theta$.
In our main Theorem we show $R_\theta$ with $|\theta_u| > 1\ \forall\ u \in U$ to be isomorphic to core regular double Stone algebras, CRDSA, that are complete and atomic, and that the crisp, {\color{Red}or definable}, sets form a complete atomistic Boolean algebra. {\color{Red}These guarantees of infimum/supremeum for arbitrary subsets and formulations in terms of fundamental elements are likely useful if dealing with equivalence relations with an infinite number of partitions, such as projective Hilbert spaces.} We further derive that every CRDSA is isomorphic to a subalgebra of a principal rough set algebra, $R_\theta$, for some approximation space $\langle U,\theta \rangle$. In our main Corollary we show explicitly how to embed $R_\theta$ into the CRDSA and first demonstrate by extending the culminating finite example of \cite{RCRDSA}. {\color{Red}As our capstone, we consider the projective Hilbert space of complex numbers, $\mathbb{C}$ and show, among other things, the power set of the set of pure states and is a complete, atomistic Boolean algebra. In closing, we suggest other Quantum relevant applications that may be useful, such as Hilbert spaces of operators.}
\end{abstract}

\maketitle
\newpage
\section{Introduction}\label{S:intro} 
In this note we prove the following result:

\begin{main} 
Given an approximation space $\langle U,\theta \rangle$ let $E$ be the indexing set for the equivalence classes of $\theta$. If $|\theta_u| > 1 \ \forall u \in U$, then its "principal rough set algebra (PRSA)" resulting from its collection of rough sets, $R_\theta$, is isomorphic to $C_3^E$, where $C_3$ is the 3-element chain, and $TP_E$, where $TP_E$ is the collection of ternary set partitions on $E$, as Core Regular Double Stone algebras. Further, the three CRDSA's in question are complete and atomic with complete atomistic Boolean centers isomorphic to $C_2^E$.
\end{main}

That result is derived from the following:

\begin{proposition}
Given an approximation space $\langle U,\theta \rangle$ let $C(R_\theta)$ denote the center of its PRSA, $R_\theta$. Assume that $E$ is the indexing set for the equivalence classes of $\theta$ and $\forall\  e \in E, y_e$ has been chosen as the representative of its equivalence class. Then $C(R_\theta) = \{(\cup\theta_{y_e},\cup\theta_{y_e})_{e \in B} : B\in \mathcal{P}(E)\}$ where $\mathcal{P}(E)$ denotes the power set of $E$. We note $B = \emptyset$ yields the rough set $(\emptyset,\emptyset)$.
\end{proposition}

\begin{proposition}
Given an approximation space $\langle U,\theta \rangle$ let $C(R_\theta)$ denote the center of its PRSA, $R_\theta$. Assume that $E$ is the indexing set for the equivalence classes of $\theta$, then $C(R_\theta) \cong 2^{E}$ as a Boolean algebra.
\end{proposition}

We give an alternate, simpler, proof of the following result from \cite{RCRDSA}:
\begin{proposition}
Given an approximation space $\langle U,\theta \rangle$, its collection of rough subsets, $R_\theta$, is a Core Regular Double Stone Algebra $\iff$ $|\theta_u| > 1 \ \forall u \in U$.
\end{proposition}

Before we move on, we make the connection to the results of \cite{dai} where they build a structure called a "3-valued Lukasiewicz algebra" which they denoted $RS(U)$ and has roots in the important Lukasiewicz logic. 

\begin{proposition}
Given an approximation space $\langle U,\theta \rangle$, its PRSA, $R_\theta$, is a Core Regular Double Stone Algebra $\iff$ it's rough 3-valued Lukasiewicz algebra, $RS(U)$, is centered, e.g. a 3-valued Post algebra.\cite{dai} Furthermore, it's center element as a 3-valued Lukasiewicz algebra is $\langle \emptyset,U \rangle$.
\end{proposition}

Then we move on to proving the following main Corollary:

\begin{mainc}
Given an approximation space $\langle U,\theta \rangle$, if $|\theta_u| > 1 \ \forall u \in U$, then its PRSA, $R_\theta$, is embedded into $C_3^U$ as a CRDSA via the following maps:\\
${\color{Red}\alpha_r \circ \phi}:R_\theta\hookrightarrow TP_U\hookrightarrow C_3^U$
where $\phi(\langle \underline{X},\overline{X} \rangle) = (\underline{X},\overline{X}^c)$ and $\alpha_r$ is $\alpha$ as defined in Theorem \ref{L:ds} restricted to the image of $\phi$, and hence we can identify it with its images.
\end{mainc}

These main results allow us to follow in the footsteps of Theorem 3. and Corollary 2.4 of \cite{comer} with the next two results:  

\begin{proposition}
Let $C_3^J$ be a CRDSA, $\langle U,\theta \rangle$ be the approximation space given by $U = J \times \{0,1\}$, $\theta = \{(j0),(j1)\} : j \in J\}$ and $R_\theta$ be the PRSA resulting from $\langle U,\theta \rangle$, then $R_\theta \cong C_3^J$.
\end{proposition} 

\begin{proposition}
Every core regular double Stone algebra is isomorphic to a subalgebra of a principal rough set algebra, $R_\theta$, for some approximation space $\langle U,\theta \rangle$.
\end{proposition}

Finally, we explicitly demonstrate the main Theorem and Corollary by expanding the culminating finite example from \cite{RCRDSA} {\color{Red}and, as our capstone, we consider the projective Hilbert space of complex numbers.}

\section{Algebras of Rough Sets}\label{CRA} 
In this section we follow the approach to rough sets followed in \cite{comer}. A pair $\langle U,\theta \rangle$ that consists of an equivalence relation $\theta$ on a nonempty set $U$ is called an approximation space. Every $X \subseteq U$ has an upper approximation $\overline{X}$ and a lower approximation $\underline{X}$ in terms of the $\theta$-classes. Namely, $\overline{X} = \bigcup \{ \theta_x : \theta_x \cap{X} \neq \emptyset \}$ and $\underline{X} = \bigcup \{ \theta_x : \theta_x \subseteq{X} \}$, a $Rough\ Set$ is the pair $\langle\underline{X},\overline{X}\rangle$. We denote $\overline{X} \setminus \underline{X} = BD(X)$ and call it the $boundary$ of $X$. We further note that $(\underline{X},BD(X),\overline{X}^c)$ is the ternary set partition of $U$ fully determined by $\langle\underline{X},\overline{X}\rangle$ and hence we can pass interchangeably between the two. From here forward we will denote $(\underline{X},BD(X),\overline{X}^c)$ as $(\underline{X},\overline{X}^c)$ since $BD(X)$ is determined by $(\underline{X},\overline{X}^c)$ in the ternary set partition representation.

Given $\langle U,\theta \rangle$, we denote the collection of all rough subsets of $U$ by $R_\theta$. In order to describe the structure placed on $R_\theta$ we need a few definitions:

\begin{definition}\label{atomistic} 
Let $P$ be a lattice with least element $0$.
\begin{enumerate}
\item $a \in P$ is an atom if $0 < a$ and $\nexists x \in P$ such that $0 < x < a$. 
\item $P$ is atomic if $\forall b \in P$, $\exists a \in P$ such that $0 < a \leq b$.
\item $P$ is atomistic if $\forall b \in P$, $b$ is the join of a set of atoms.
\end{enumerate}
\end{definition}

\begin{definition}\label{DSA}A double Stone algebra (DSA) $<L,\wedge,\vee,*,+,0,1>$ is an algebra of type $<2,2,1,1,0,0>$ such that:
\begin{enumerate}
	\item $<L,\vee,\wedge,0,1>$ is a bounded distributive lattice
	\item $x^*$ is the pseudocomplement of $x$ i.e $y \leq x^* \leftrightarrow y \wedge x = 0$
	\item $x^+$ is the dual pseudocomplement of $x$ i.e $y \geq x^+ \leftrightarrow y \vee x = 1$
	\item $x^* \vee x^{**}=1$, $x^+ \wedge x^{++}=0$, e.g. the Stone identities.
\end{enumerate}
\end{definition}
Condition 4. is what distinguishes from double p-algebra and conditions 2. and 3. are equivalent to the equations
\begin{itemize}
	\item $x \wedge (x \wedge y)^* = x \wedge y^*$, $x \vee (x \vee y)^+ = x \vee y^+$
	\item $x \wedge 0^* = x$, $x \vee 1^+ = x$
	\item $0^{**} = 0$, $1^{++} = 1$
\end{itemize}
so that DSA is an equational class.

\begin{definition}\label{RDSA} A double Stone algebra L is called regular, if it additionally satisfies
\begin{itemize}
	\item $x \wedge x^+ \leq y \vee y^*$
	\begin{itemize}
		\item this is equivalent to $x^+ = y^+$ and $x^* = y^* \rightarrow x = y$
	\end{itemize}
\end{itemize}
We will use the abbreviation RDSA to denote regular double Stone Algebra 
\end{definition}
\begin{note}
Regular means any 2 congruences having a common class are the same.
\end{note}

\begin{definition}\label{PRSA} Given an approximation space $\langle U,\theta \rangle$ and it's collection of rough subsets $R_\theta$, we define the "principal rough set algebra" it determines, $R_\theta = \langle R_\theta,\vee,\wedge,*,+,0,1\rangle$, as follows:
\begin{itemize}
	\item $0 = \langle \emptyset,\emptyset \rangle$
	\item $1 = \langle U,U \rangle$
	\item $\langle\underline{X},\overline{X}\rangle \vee \langle\underline{Y},\overline{Y}\rangle = \langle\underline{X} \cup \underline{Y},\ \overline{X} \cup \overline{Y}\rangle$
	\item $\langle\underline{X},\overline{X}\rangle \wedge \langle\underline{Y},\overline{Y}\rangle = \langle\underline{X} \cap \underline{Y},\ \overline{X} \cap \overline{Y}\rangle$
	\item $\langle\underline{X},\overline{X}\rangle^* = \overline{X}^c,\overline{X}^c \rangle$
	\item $\langle\underline{X},\overline{X}\rangle^+ = \langle \underline{X}^c,\underline{X}^c \rangle$
\end{itemize}
It was Pomykala in \cite{pomykala} who first discovered that under $\vee, \ \wedge$ and $^* \ R_\theta$ forms a Stone algebra. This was improved upon by Comer in \cite{comer} by noting that with the addition of $^+ \ R_\theta$ actually forms a regular double Stone algebra. Finally, we note that Dunstch \cite{duns} recognized this algebra followed a special form that occurred more generally which he referred to as Katrinak algebras.

For the purposes of this note we will refer to this specific structure formed out of a collection of rough subsets $R_\theta$ as the \textit{principal rough set algera (PRSA)} determined by $\langle U,\theta \rangle$ and will use $R_\theta$ to denote it and the corresponding collection of rough sets interchangeably.    
\end{definition}

\begin{definition}\label{center}
Let L be a regular double Stone algebra. An element $x \in L$ is called central if $x^*=x^+$. The set of all central elements of L is called the center of L and is denoted by C(L); that is $C(L)=\{x \in L|x^*=x^+\}=\{x^* |x \in L\}=\{x^+ |x \in L\}$.  
\end{definition}

Now we cite the following Theorems and Definitions from \cite{CRDSA}: %\textit{Centre of Core Regular Double Stone Algebra} \cite{CRDSA}:

\begin{theorem}\label{centerboolsub} 
Let L be a regular double Stone algebra. Then $C(L)$ is a Boolean sub algebra of L with respect to the induced operations $\wedge$, $\vee$ and $^*$.\cite{CRDSA}
\end{theorem}

\begin{definition}\label{dense}
Every element x of a double Stone algebra L with the property $x^* = 0$ (or equivalently, $x** =1$) is called dense. Every element of the form $x \vee x^*$ is dense and we denote the set of all dense elements of L, D(L). Every element x with the property $x^+ = 1$ (or equivalently, $x^{++} = 0$) is called dually dense. Every element of the form $x \wedge x^+$ is dually dense and we denote the set of all dual dense elements of L, $\overline{D(L)}$.
\end{definition}

\begin{definition}\label{core}
The core of a double Stone algebra L is defined to be $K(L) = D(L) \cap \overline{D(L)}$ and we call a regular double Stone algebra with non-empty core a core regular double Stone algebra, denoted CRDSA.\cite{CRDSA}
\end{definition}

For any CRDSA L, $|K(L)| = 1$ follows easily from regularity and we will call this element $h$. We now cite 2 useful results from \cite{CRDSA} that we will use in the next section.

\begin{theorem}\label{coreexp}
If A is CRDSA with core element h, then every element x of A can be written as $x = x^{**} \wedge (x^{++} \vee h)$
and $x = x^{++} \vee (x^{**} \wedge h)$.
\end{theorem}

\begin{theorem}\label{centiso}\cite{CRDSA}
Two CRDSAs are isomorphic if and only if their centers are isomorphic.
\end{theorem}

It is the CRDSA structure that we seek to leverage in further refining the structure of $R_\theta$ and we define CRDSA as an algebra of type $<2,2,1,1,0,0,0>$ where the core element $h$ is also a constant.  

\begin{definition}\label{CoreCon}
A core regular double Stone algebra is an algebra $<L,\wedge,\vee,*,+,0,1,h>$ of type $<2,2,1,1,0,0,0>$ where $<L,\wedge,\vee,*,+,0,1>$ is a regular double Stone algebra and $h$ is its unique core element.
\end{definition}

\begin{definition}\label{NSBDLdef}
   Let $J$ be a non-empty set and let $L = \{(X_1,X_2) : X_1,X_2 \subseteq J$ and $X_1 \cap X_2 = \emptyset \}$. We define binary operations $\vee$ and $\wedge$ on L as follows:
\begin{itemize}
	\item $(X_1,X_2) \vee (Y_1,Y_2) = (X_1 \cup Y_1,X_2 \cap Y_2)$ and \label{m-i}
	\item $(X_1,X_2) \wedge (Y_1,Y_2) = (X_1 \cap Y_1,X_2 \cup Y_2)$ \label{j-i}
	\begin{itemize}
		\item The fact that $L$ is a bounded distributive lattice with bounds $(J,\emptyset)$ and $(\emptyset,J)$ is well known.\label{c-s}
	\end{itemize}
   \end{itemize}
We call this the "ternary partition lattice" and denote it $TP_J$. We note that as $(X_1,X_2)$ determine $(X_1 \cup X_2)^c$, we may not list it and that in addition to $(J,\emptyset)$ and $(\emptyset,J)$, there is another very important element of $TP_J$, namely $h=(\emptyset,\emptyset)$. Lastly, we note that $\vee$ and $\wedge$ coincide with the following partial ordering on L:
\begin{itemize}
	\item $(X_1,X_2) \leq (Y_1,Y_2) \leftrightarrow X_1 \subseteq Y_1$ and $Y_2 \subseteq X_2$
\end{itemize}
\end{definition}

Let $C_3$ denote the 3-element chain with every element a constant. First we need to mention some facts about ternary set partitions and the 3 element chain $C_3$.

\begin{definition}\label{C3alg}
Let $C_3\ = \ \{0,h,1\}$ be the 3 element chain and define it to be an algebra $<C_3,\wedge,\vee,*,+,0,1,h>$ of type $<2,2,1,1,0,0,0>$ where the operations $\vee$, $\wedge$ are the lattice join, meet operations, respectively and all 3 elements are constants. Further define $*,+$ as follows:
\begin{itemize}
	\item $* : C_3 \rightarrow C_3$ is defined by $0 \rightarrow 1$, $h \rightarrow 0$ and $1 \rightarrow 0$ and
	\item $+ : C_3 \rightarrow C_3$ is defined by $0 \rightarrow 1$, $h \rightarrow 1$ and $1 \rightarrow 0$
\end{itemize} 
\end{definition}

First we note that $C_3$ and hence $C_3^J$ are CRDSA, ref.\cite{gamesec} and from here forward we treat it as such. We finish this section by including 2 helpful figures and by citing some useful representation results we will leverage in the next section.\cite{gamesec},\cite{nearlybool} 

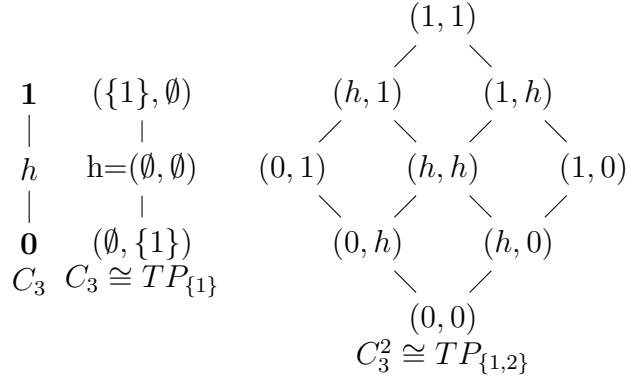
\begin{figure}\label{c3}
\caption{Minimal Examples, $C_3$ and $C_3^2$}
    \begin{center}
    \begin{tikzpicture}
        \tikzstyle{every node} = [rectangle]
	%\node (s0) at (0,.5) {$C_3^2\cong TP^{\{1,2\}}$};
        %\node (s1) at (0,1) {($\emptyset,\{1,2\}$)};
        %\node (s2) at (-1,2) {($\emptyset,\{1\}$)};
        %\node (s3) at (1,2) {($\emptyset,\{2\}$)};
        %\node (s4) at (-2,3) {($\{2\},\{1\}$)};
	%\node (s5) at (0,3) {h=($\emptyset,\emptyset$)};
	%\node (s6) at (2,3) {($\{1\},\{2\}$)};     
        %\node (s7) at (-1,4) {($\{2\},\emptyset$)};
	%\node (s8) at (1,4) {($\{1\},\emptyset$)};
	%\node (s9) at (0,5) {($\{1,2\},\emptyset$)};
	\node (s0) at (0,.5) {$C_3^2\cong TP_{\{1,2\}}$};
        \node (s1) at (0,1) {($0,0$)};
        \node (s2) at (-1,2) {($0,h$)};
        \node (s3) at (1,2) {($h,0$)};
        \node (s4) at (-2,3) {($0,1$)};
	\node (s5) at (0,3) {($h,h$)};
	\node (s6) at (2,3) {($1,0$)};     
        \node (s7) at (-1,4) {($h,1$)};
	\node (s8) at (1,4) {($1,h$)};
	\node (s9) at (0,5) {($1,1$)};
	%\node (s18) at (5.5,.5) {$C_3^2$};
	%\node (s19) at (5.5,1) {($0,0$)};
       %\node (s20) at (4.5,2) {($0,h$)};
        %\node (s21) at (6.5,2) {($h,0$)};
        %\node (s22) at (3.5,3) {($0,1$)};
	%\node (s23) at (5.5,3) {($h,h$)};
	%\node (s24) at (7.5,3) {($1,0$)};     
        %\node (s25) at (4.5,4) {($h,1$)};
	%\node (s26) at (6.5,4) {($1,h$)};
	%\node (s27) at (5,5) {($1,1$)};
	\node (s13) at (-4,1.5) {$C_3\cong TP_{\{1\}}$};
	\node (s10) at (-4,2) {($\emptyset,\{1\}$)};
	\node (s11) at (-4,3) {h=($\emptyset,\emptyset$)};
	\node (s12) at (-4,4) {($\{1\},\emptyset$)};
	\node (s17) at (-5.5,1.5) {$C_3$};
	\node (s14) at (-5.5,2) {$\mathbf{0}$};
	\node (s15) at (-5.5,3) {$h$};
	\node (s16) at (-5.5,4) {$\mathbf{1}$};
        \foreach \from/\to in {s1/s2, s1/s3, s2/s4, s2/s5, s3/s5, s3/s6, s4/s7, s5/s7, s5/s8, s6/s8, s7/s9, s8/s9, s10/s11, s11/s12, s14/s15, s15/s16}%, s19/s20, s19/s21, s20/s22, s20/s23, s21/s24, s21/s23, s22/s25, s23/s25, s24/s26, s23/s26, s25/s27, s26/s27}
            \draw[-] (\from) -- (\to);
    \end{tikzpicture}
    \end{center}
\end{figure}

\begin{theorem}\label{L:ds}
$TP_J \cong C_3^J$ as CRDSA under the following mapping:\cite{gamesec}\\
$\alpha : TP_J \rightarrow C_3^J$ pointwise for $(X_1,X_2) \in TP_J$ by $\alpha(X_1,X_2)=(x_i)_{i \in J}=$
\begin{enumerate}
	\item 1 if $i \in X_1$
	\item h if $i \in (X_1 \cup X_2)^c$
	\item 0 if $i \in X_2$
\end{enumerate}
\end{theorem}

Now recall that the ternary set partion $(\underline{X},\overline{X}^c)$ of $U$ is fully determined by $\langle\underline{X},\overline{X}\rangle$ and hence we will be looking to leverage $\alpha$ to represent $\langle\underline{X},\overline{X}\rangle$ as it's corresponding element of $C_3^U$ in such a way that the CRDSA structure is preserved. We further note that every CRDSA is a subdirect product of $C_3$, similarly as for Boolean algebras and $C_2$, and that the variety generated by $C_3$ is dually equivalent to the category of Stone spaces and hence the category of Boolean algebras.\cite{nearlybool}

\section{Some New and Some Old}\label{somenew}

\begin{theorem}\label{CPRSA}
Given an approximation space $\langle U,\theta \rangle$, its PRSA, $R_\theta$, is a Core Regular Double Stone Algebra $\iff$ $|\theta_u| > 1 \ \forall u \in U$.
\end{theorem}
We note that most of the work has been done for us as we already have that $R_\theta$ is a RDSA.\cite{pomykala}, \cite{comer}, see Definition \ref{PRSA} discussion. 
\begin{proof}
$\Longrightarrow$ If we assume $\exists$ a set $X \subseteq U$ such that $(\underline{X},\overline{X})=(\emptyset,U)$, it follows from the definition of $(\underline{X},\overline{X})$ that $|\theta_u|\ >\ 1\  \forall\ u \in U$.\\
$\Longleftarrow$ Assume $\forall\ u \in U, \ |\theta_u|\ >\ 1$, that $E$ is the indexing set for the equivalence classes of $\theta$ and $\forall\  e \in E, y_e \in U$ has been chosen as the representative of its equivalence class. We will construct a set $X$ such that $(\underline{X},\overline{X})=(\emptyset,U)$, e.g. is both dense and dually dense. $\forall \theta_{y_e}$, choose $x_e \in \theta_{y_e}$ with $x_e \neq y_e$ and let $X = \cup \{x_e\}_{e \in E}$. By construction, $(\underline{X},\overline{X})=(\emptyset,U)$.\\ 
Note: When $R_\theta$ is a CRDSA we consider it to be the algebra\\
$<R_\theta,\wedge,\vee,*,+,0,1,h>$ of type $<2,2,1,1,0,0,0>$
\end{proof}

We further note that in the first part of the above proof there could be many sets $X$ that will yield the core element of $R_\theta$, but all will yield the ternary set representation is $(\emptyset,U)$ and its image under $\alpha$ from Theorem \ref{L:ds} is $(h)_{i \in U}$. Before we move on to our Main Theorem we wish to list a few more results that follow directly from the results of Theorem \ref{CPRSA} and the results of \cite{dai}, but using the perspective developed here. There, given an approximation space $\langle U,\theta \rangle$ many nice results are derived by placing a different structure on $R_\theta$. They, and others before them, called that structure a "3-valued Lukasiewicz algebra" which they denoted $RS(U)$, has roots in the important Lukasiewicz logic and hence we will note a the connection with CRDSA here. We will not go through the details of the structure other than to say it has different unary operators and refer the interested reader to \cite{dai} for more information.

\begin{corollary}\label{3VLA}
Given an approximation space $\langle U,\theta \rangle$, its PRSA, $R_\theta$, is a Core Regular Double Stone Algebra $\iff$ it's rough 3-valued Lukasiewicz algebra, $RS(U)$, is centered, e.g. a 3-valued Post algebra.\cite{dai} Furthermore, it's center element as a 3-valued Lukasiewicz algebra is $\langle \emptyset,U \rangle$.
\end{corollary}
\begin{proof} 
This follows immediately from Theorem \ref{CPRSA} above and Theorem 5 of \cite{dai}.
\end{proof}

Assume that $E$ is the indexing set for the equivalence classes of $\theta$, we wish to show $R_\theta \cong C_3^E$ as CRDSA. We will first show that if $X \subseteq U$ is such that $(\underline{X},\overline{X}) \in C(R_\theta)$ then $X = \underline{X} = \overline{X}$ and characterize $C(R_\theta)$ as a Boolean algebra. Given $R_\theta$, we note that sets $X \subseteq U$ such that  $X = \underline{X} = \overline{X}$ are classically referred to in Rough Set theory literature as crisp or definable sets.

\begin{lemma}\label{PRSACenter}
Given an approximation space $\langle U,\theta \rangle$ let $C(R_\theta)$ denote the center of its PRSA, $R_\theta$. Assume that $E$ is the indexing set for the equivalence classes of $\theta$ and $\forall\  e \in E, y_e \in U$ has been chosen as the representative of its equivalence class. Then $C(R_\theta) = \{\langle \cup\theta_{y_e},\cup\theta_{y_e} \rangle_{e \in B} : B\in \mathcal{P}(E)\}$ where $\mathcal{P}(E)$ denotes the power set of $E$. We note $B = \emptyset$ yields the rough set $(\emptyset,\emptyset)$.
\end{lemma}
\begin{proof} 
$\Longrightarrow$ If we assume $\langle \underline{X},\overline{X})\in C(R_\theta \rangle$, then $\overline{X} = \cup\{ \theta_x : \theta_x \cap{X} \neq \emptyset \} = \underline{X}= \cup\{ \theta_x : \theta_x \subseteq{X} \}= X$ and $\theta_{y_e}\cap X \neq \emptyset,\ \Longrightarrow\ \theta_{y_e} \subseteq X\ \Longrightarrow \underline{X}=\cup\{\theta_{y_e} : \theta_{y_e}\cap X \neq \emptyset\}$.\\
$\Longleftarrow$ $(\underline{X},\overline{X}) = \langle \cup\theta_{y_e},\cup\theta_{y_e}\rangle_{e \in B}$ for some $B\in \mathcal{P}(E)\Longrightarrow (\underline{X},\overline{X}) \in C(R_\theta)$.
\end{proof}

We note that this result yields the characterization of $C(R_\theta)$ as a Boolean algebra and as the subset of crisp sets of $R_\theta$.

\begin{corollary}\label{PRSACenterBool}
Given an approximation space $\langle U,\theta \rangle$ let $C(R_\theta)$ denote the center of its PRSA, $R_\theta$. Assume that $E$ is the indexing set for the equivalence classes of $\theta$, then $C(R_\theta)$ is exactly the set of crisp sets of $R_\theta$ and is complete, atomistic and isomorphic to $2^{E}$ as a Boolean algebra.
\end{corollary}
\begin{proof}
The result follows immediately from Theorem \ref{centerboolsub}, Lemma \ref{PRSACenter} and the fact that the power set is a complete atomistic Boolean algebra. 
\end{proof}

\begin{theorem}\label{C3CPRSA}\textbf{The Main Theorem}\\
Given an approximation space $\langle U,\theta \rangle$ let $E$ be the indexing set for the equivalence classes of $\theta$. If $|\theta_u| > 1 \ \forall u \in U$, then its PRSA, $R_\theta$, is isomorphic to $C_3^E$, where $C_3$ is the 3-element chain, and $TP_E$, where $TP_E$ is the collection of ternary set partitions on $E$, as Core Regular Double Stone algebras. Further, the three CRDSA's in question are complete and atomic with complete atomistic Boolean centers isomorphic to $C_2^E$.  
\end{theorem}
\begin{proof}
Corollary \ref{PRSACenterBool} characterizes $C(R_\theta) \cong 2^{E}$ as a complete atomistic Boolean algebra, Theorem \ref{centiso} establishes that $R_\theta \cong C_3^E$ and Theorem \ref{L:ds} establishes the final isomorphism with $TP_E$.  $C_3^E$ is clearly atomic and next we show $TP_E$ is complete.\\ 
Let $D$ be an indexing set for the pairs $\{(X_d,Y_d) : \{X_d,Y_d\} \subseteq P(E), X_d \cap Y_d = \emptyset\}$, $I \subseteq D$, $Z_I = \{(X_i,Y_i)_{i \in I}\} \subseteq TP_E$, $J_I = \cup(X_i)_{i \in I}$ and $M_I = \cap(Y_i)_{i \in I}$. Now we consider the cases for $\wedge$ and $\vee$ from Definition \ref{NSBDLdef}.
\begin{enumerate}
	\item $(X_1,Y_1) \vee (X_2,Y_2) = (X_1 \cup X_2,Y_1 \cap Y_2)$, consider $\vee Z_I$
		\begin{itemize}
			\item $J_I$ is the least upper bound of $((X_i)_{i \in I},\subseteq)$, $M_I$ is the greatest lower bound of $((Y_i)_{i \in I},\subseteq)$ and hence $\vee Z_I = (J_I,M_I) = (\cup(X_i)_{i \in I},\cap(Y_i)_{i \in I})$
		\end{itemize}
		\item $(X_1,Y_1) \vee (X_2,Y_2) = (X_1 \cup X_2,Y_1 \cap Y_2)$, similarly $\wedge Z_I = (\cap(X_i)_{i \in I},\cup(Y_i)_{i \in I})$
\end{enumerate}
 %as $\vee$, $\wedge$ involve only $\cup$ and $\cap$, for any $S \subseteq E$, $\cap S$ is the infimum, $\cup S$ is the supremum of $(S,\subseteq)$ and $(X_1,X_2) \leq (Y_1,Y_2 \leftrightarrow X_1 \subseteq Y_1$ and $Y_2 \subseteq X_2$ determines the partial order of $TP_E$.$%
\end{proof}

In light of the Main Theorem we wish to note a couple more results that follow easily. Recall that in \cite{nearlybool} we noted that every CRDSA is a subdirect product of $C_3^J$ for some indexing set $J$ as a result of Theorem 2. of \cite{comer}. We now unabashedly now follow in the footsteps of Theorem 3. and Corollary 2.4 of \cite{comer}.  

\begin{theorem}\label{CRDSAtoPRSA}
Let $C_3^J$ be a CRDSA, $\langle U,\theta \rangle$ be the approximation space given by $U = J \times \{0,1\}$, $\theta = \{(j0),(j1)\} : j \in J\}$ and $R_\theta$ be the PRSA resulting from $\langle U,\theta \rangle$, then $R_\theta \cong C_3^J \cong TP_J$ and all 3 CRDSA are complete and atomic.
\end{theorem} 
\begin{proof}
One only need apply the Main Theorem.
\end{proof}

\begin{corollary}
Every core regular double Stone algebra is isomorphic to a subalgebra of a principal rough set algebra, $R_\theta$, for some approximation space $\langle U,\theta \rangle$.
\end{corollary}
\begin{proof}
This follows from Theorem \ref{CRDSAtoPRSA} and the fact that every CRDSA is a subdirect product of $C_3$.
\end{proof}

We note that \cite{RCRDSA} and \cite{nearlybool} contain some useful results for finite $E$. Now we will use the ternary set partiton representation of $R_\theta$ to embed it into $C_3^U$ as a CRDSA.

\begin{corollary}\label{CPRSAinC3U}\textbf{The Main Corollary}\\
Given an approximation space $\langle U,\theta \rangle$, if $|\theta_u| > 1 \ \forall u \in U$, then its PRSA, $R_\theta$, is embedded into $C_3^U$ as a CRDSA via the following maps:\\
${\color{Red}\alpha_r \circ \phi}:R_\theta\hookrightarrow TP_U\hookrightarrow C_3^U$
where $\phi(\langle \underline{X},\overline{X} \rangle) = (\underline{X},\overline{X}^c)$ and $\alpha_r$ is $\alpha$ as defined in Theorem \ref{L:ds} restricted to the image of $\phi$.
\end{corollary}
\begin{proof}
Theorem \ref{CPRSA} establishes $R_\theta$ as a CRDSA and as $\alpha:TP_U\longrightarrow C_3^U$ is a CRDSA isomorphism by Theorem \ref{L:ds}, all we need to show is that $\phi$ is an injective CRDSA homomorphism.
\begin{itemize}
	\item $\phi(\langle \emptyset,\emptyset \rangle) = (\emptyset,U)$, hence $\phi$ preserves $0$
	\item $\phi(\langle U,U \rangle) = (U,\emptyset)$, hence $\phi$ preserves $1$
	\item $\phi(\langle \emptyset,U \rangle) = (\emptyset,\emptyset)$, hence $\phi$ preserves $h$
	\item $\phi(\langle \underline{X},\overline{X} \rangle^*) = \phi(\langle \overline{X}^c,\overline{X}^c \rangle) = (\overline{X}^c,\overline{X}) = (\underline{X},\overline{X}^c)^*$, hence $\phi$ preserves $^*$
	\item $\phi(\langle \underline{X},\overline{X} \rangle^+) = \phi(\langle \underline{X}^c,\underline{X}^c \rangle) = (\underline{X}^c,\underline{X}) = (\underline{X},\overline{X}^c)^+$, hence $\phi$ preserves $^+$
	\item $\phi(\langle \underline{X},\overline{X} \rangle \vee \langle\underline{Y},\overline{Y} \rangle) = \phi(\langle\underline{X}\cup\underline{Y},\overline{X}\cup\overline{Y} \rangle) = (\underline{X}\cup\underline{Y},(\overline{X}\cup\overline{Y})^c) =\\ (\underline{X}\cup\underline{Y},\overline{X}^c\cap\overline{Y}^c) = (\underline{X},\overline{X}^c)\vee(\underline{Y},\overline{Y}^c)$, hence $\phi$ preserves $\vee$
	\item $\phi(\langle \underline{X},\overline{X} \rangle \wedge \langle \underline{Y},\overline{Y} \rangle) = \phi(\langle \underline{X}\cap\underline{Y},\overline{X}\cap\overline{Y} \rangle) = (\underline{X}\cap\underline{Y},(\overline{X}\cap\overline{Y})^c) =\\ (\underline{X}\cap\underline{Y},\overline{X}^c\cup\overline{Y}^c) = (\underline{X},\overline{X}^c)\wedge(\underline{Y},\overline{Y}^c)$, hence $\phi$ preserves $\wedge$
	\item $(\underline{X},\overline{X}^c) = (\underline{Y},\overline{Y}^c) \iff (\underline{X},\overline{X}) = (\underline{Y},\overline{Y})$, hence $\phi$ is injective.
\end{itemize} 
\end{proof}

\begin{note}
Perhaps the Main Corollary is somewhat expected, in light of the fact that every CRDSA is a subdirect product of $C_3$.\cite{nearlybool} However, from it we can specifically "trace" $\langle \underline{X},\overline{X} \rangle \hookrightarrow (\underline{X},\overline{X}^c) \in TP_U \hookrightarrow \alpha((\underline{X},\overline{X}^c)) \in C_3^U$ and in the next section we demonstrate the Main Corollary and Main Theorem. For completeness we also note that the approximation space from Theorem \ref{CRDSAtoPRSA} will give us back the approximation space we start with in the following example.  
\end{note}

\section{Examples of the Main Theorem and Corollary}
As a first example that reflects the note ending the last section we would like to continue the analysis of example 3.6, $A_9$, from \cite{RCRDSA}, which is $C_3^2$ from our Table 1. We encourage readers to read \cite{RCRDSA} in it's entirety but here we will build on the initial analysis from pages 338 and 339 and expand on their Table 1. In order to reflect the perspectives of the Main Theorem and Main Corollary of this note we will only need consider a universe of cardinality 4 and an equivalence relation that partitions the universe into 2 subsets of cardinality 2 for the example given on the next page.

\begin{example} 
Consider the universe $X = \{P_a,P_1a,P_b,P_1b\}$, equivalence relation $\sim = \{\{P_a,P1a\},\{P_b,P_1b\}\}$ from pages 338 and 339 of \cite{RCRDSA}. We wish to abstract away from that to the perspectives of the Main Theorem and the Main Corollary. Hence we will identify $X$, $\sim$ with $\langle U, \theta \rangle$ given by $U = \{w,x,y,z\}$ and $\theta = \{\{w,x\},\{y,z\}\}$ and consider the PRSA, ${R_\theta}$, from Definition \ref{PRSA}. We now will create the Table $1_{R_\theta}$ which is an extension of Table 1 on page 338 of \cite{RCRDSA} that is formed by adding the columns necessary to demonstrate the Main Theorem and the Main Corollary.

\begin{table}[h!]
\centering
\caption{Table $1_{R_\theta}$, an Example of the Main Theorem and the Main Corollary}
%\vspace{3mm}
\begin{tabular}{||c c c c c c c||} 
 \hline
 $X\subseteq U$ & $\underline{X}$ & $\overline{X}$ & $\overline{X}^c$ & $TP_U$ & ${C_3}^U$ & ${C_3}^2$ \\ [0.5ex] 
 \hline\hline
 $\emptyset$ & $\emptyset$ & $\emptyset$ & X & $(\emptyset,X)$ & (0000) & (00) \\ 
 w & $\emptyset$ & $\{w,x\}$ & $\{y,z\}$ & $(\emptyset,\{y,z\})$ & (hh00) & (h0) \\
 x & $\emptyset$ & $\{w,x\}$ & $\{y,z\}$ & $(\emptyset,\{y,z\})$ & (hh00) & (h0) \\
 y & $\emptyset$ & $\{y,z\}$ & $\{w,x\}$ & $(\emptyset,\{w,x\})$ & (00hh) & (0h) \\
 z & $\emptyset$ & $\{y,z\}$ & $\{w,x\}$ & $(\emptyset,\{w,x\})$ & (00hh) & (0h) \\
 $\{w,y\}$ & $\emptyset$ & X & $\emptyset$ & $(\emptyset,\emptyset)$ & (hhhh) & (hh) \\ 
 $\{w,x\}$ & $\{w,x\}$ & $\{w,x\}$ & $\{y,z\}$ & $(\{w,x\},\{y,z\})$ & (1100) & (10) \\
 $\{w,z\}$ & $\emptyset$ & X & $\emptyset$ & $(\emptyset,\emptyset)$ & (hhhh) & (hh) \\
 $\{y,x\}$ & $\emptyset$ & X & $\emptyset$ & $(\emptyset,\emptyset)$ & (hhhh) & (hh) \\
 $\{y,z\}$ & $\{y,z\}$ & $\{y,z\}$ & $\{w,x\}$ & $(\{y,z\},\{w,x\})$ & (0011) & (01) \\ 
 $\{x,z\}$ & $\emptyset$ & X & $\emptyset$ & $(\emptyset,\emptyset)$ & (hhhh) & (hh) \\ 
 $\{w,x,z\}$ & $\{w,x\}$ & X & $\emptyset$ & $(\{w,x\},\emptyset)$ & (11hh) & (1h) \\
 $\{w,x,y\}$ & $\{w,x\}$ & X & $\emptyset$ & $(\{w,x\},\emptyset)$ & (11hh) & (1h) \\
 $\{w,z,y\}$ & $\{y,z\}$ & X & $\emptyset$ & $(\{y,z\},\emptyset)$ & (hh11) & (h1) \\
 $\{z,x,y\}$ & $\{y,z\}$ & X & $\emptyset$ & $(\{y,z\},\emptyset)$ & (hh11) & (h1) \\
 X & X & X & $\emptyset$ & $(X,\emptyset)$ & (1111) & (11) \\ [1ex] 
 \hline
\end{tabular}
%\vspace{3mm}
%\caption{Table to test captions and labels}
\label{tableRtheta}
\end{table}

In Table $1_{R_\theta}$ one can "trace" $\langle \underline{X},\overline{X} \rangle \hookrightarrow (\underline{X},\overline{X}^c) \in TP_U \hookrightarrow \alpha((\underline{X},\overline{X}^c)) \in C_3^U$ and clearly see the $R_\theta \hookrightarrow TP_U \hookrightarrow {C_3}^U$ as ${C_3}^2$. Lastly, we note that the approximation space from Theorem \ref{CRDSAtoPRSA} will give us back the approximation space we start with in the example. $|J|=2$ so let $j=\{j_0,j_1\}$, $U = \{(j_0,0),(j_0,1),(j_1,0),(j_1,1)\}$ and let $\theta = \{\{(j_0,0),(j_0,1)\},\{(j_1,0),(j_1,1)\}\}$ and $C_3^2 \cong R_\theta$ as shown in Table $1_{R_\theta}$. 

\end{example}

{\color{Red}As our culminating example, we consider the projective Hilbert space over the complex numbers $\mathbb{C}$ and the concept of the states from quantum mechanics. We only take the definitions from the quantum concept of a physical system that are necessary for us to apply our results, for more details we refer the interested reader to \cite{quant}\cite{quantc}.

\begin{definition}\label{projhilb}
Projective Hilbert Space
\begin{enumerate}
\item A Hilbert space is a vector space $H$ with an inner product $<x,y>$ such that $H$ is a complete metric space under the norm defined by $\sqrt{<x,x>}$.
\item The projective Hilbert space $P(H)$ of $H$ is the set of equivalence classes of non-zero vectors $x \in H^*$ for the equivalence relation $\theta$ given by $x \theta y \iff x=cy$ for some $c \in H^*$. It is customary to choose $x$ with $||{x}||=1$ as class representatives for $[x]$.
\end{enumerate}
It is well known that $\mathbb{C}^*$ with $<x,y>$ the vector dot product of $x$ and the complex conjugate of $y$ and the equivalence relation from Definition \ref{projhilb} forms a projective Hilbert space. The equivalence classes of $P(H)$ are referred to as rays or pure states in quantum physics.\cite{quant}
\end{definition}
}

%NOTE: I had to break the coloring statements up as latex had issues coloring across everything

\begin{example}
{\color{Red}We consider $\langle \mathbb{C},\theta\rangle$ as an approximation space. Then for any $X \subseteq \mathbb{C}^*$ we have the following:
\begin{enumerate}
\item $\underline{X} = \cup\{x \in X|[x] \subseteq X\}$, $\overline{X} = \cup\{[x] | x \in X\}$ and hence 
\item $\overline{X}^c = \cup\{[z] | z \notin X\}$ and $\phi$ yields $(\underline{X}, \cup\{[z] | z \notin X\})$  in $TP_{\mathbb{C}^*}$.
\item We note that $BD(X) = (\underline{X} \cup  \overline{X}^c)^c = \cup\{x \in X|[x] \nsubseteq X\}$. 
\end{enumerate}
We see that $BD(X)$ is the union of all pure states that intersect non-trivially with $X$, but are not contained in it. All crisp sets are unions of pure states and conversely. Furthermore, we know by the results shown here that the center of $TP_{\mathbb{C}^*}$ is exactly the power set of the set of pure states and is a complete, atomistic Boolean algebra. If we consider $\alpha$ as defined in Lemma \ref{L:ds}, we see that any pure states that are in $BD(X)$ will be pointwise mapped to element $h$ of $C_3$. We find these results very interesting and as sets of states are integral to Quantum channels, we hope results like this may be useful to the Quantum community. We further note that in Quantum the space of operators can also be Hilbert spaces. Lastly, and perhaps importantly, we note that in \cite{quantc}, many references are made to \textit{the statistical duality of states and effects}. We are not versed in said duality, but we note the following results stated there:\cite{quantc}
\begin{enumerate}
\item There is a subset $L \subset E$ of effects, called propositions or sharp effects, which has a structure of a partially ordered, orthocomplemented, orthomodular, complete lattice, with the universal bounds 0 and 1, which is atomistic, separable, has the covering property and is irreducible.
\item The set S of states can be viewed as a $\sigma$-convex set of probability measures on L, which has a sufficient set $ex(S)$ of pure states: for all $a,b \in L$, $a \leq b$ if $\alpha(a) \leq \alpha(b)$ for all $\alpha \in ex(S)$.
\item There is a bijective correspondence between the sets $ex(S)$, the pure states of S, and At(L), the atoms of L, given by the support projection $\alpha \rightarrow s(\alpha)$ with $s(\alpha)$ being the smallest element for which $\alpha(b)=1$, $b \in L$.
\end{enumerate}

Given the above, we feel it is likely that there is more to say about the results we have touched upon here and Quantum mechanics.  We welcome any discussions the Quantum community may desire, but for now we will leave further investigation into these topics as potential future research.}
\end{example}

%\newpage
\section{Conclusion}
\subsection{Mathematical Results}
Given an approximation space $\langle U,\theta \rangle$, assume that $E$ is the indexing set for the equivalence classes of $\theta$ and let $R_\theta$ denote the collection of rough sets of the form $\langle\underline{X},\overline{X}\rangle$ as a regular double Stone algebr.\cite{comer} We give an alternate, simpler, proof from the one given in \cite{RCRDSA} of the fact that if $|\theta_u| > 1\ \forall\ u \in U$ then $R_\theta$ is a core regular double Stone algebra, CRDSA. Further let $C_3$ denote the 3 element chain as a CRDSA and $TP_U$ denote the collection of ternary partitions over the set $U$. In our main Theorem we show $R_\theta$ with $|\theta_u| > 1\ \forall\ u \in U$ to be isomorphic to $TP_E$ and $C_3^E$ and that the three CRDSA's in question are complete and atomic. Further, we observe that the subset of crisp sets of $R_\theta$ form a the complete atomistic Boolean algebra isomorphic to $C_2^E$, {\color{Red}likely quite useful when dealing with a specific $R_\theta$ with infinite partitions.} We also establish the relationship between such CRDSA and their rough 3-valued Lukasiewicz algebra counterparts, e.g. they're 3-valued Post algebras.\cite{dai} We then derive that every core regular double Stone algebra is isomorphic to a subalgebra of a principal rough set algebra, $R_\theta$, for some approximation space $\langle U,\theta \rangle$. In our main Corollary we show explicitly how we can embed such $R_\theta$ in $TP_U$, $C_3^U$, respectively, ${\color{Red}\alpha_r \circ \phi}:R_\theta\hookrightarrow TP_U\hookrightarrow C_3^U$, and hence identify it with its specific images. Following in the footsteps of \cite{comer}, we show $C_3^J \cong R_\theta$ for $\langle U,\theta \rangle$ the approximation space given by $U = J \times \{0,1\}$, $\theta = \{(j0),(j1)\} : j \in J\}$ and every CRDSA is isomorphic to a subalgebra of a principal rough set algebra, $R_\theta$, for some approximation space $\langle U,\theta \rangle$. Further, we know subalgebras of $TP_E$ and $C_3^E$ must exist for every $E$ that is an indexing set for the equivalence classes of any equivalence relation $\theta$ on $U$ satisfying $|\theta_u| > 1\ \forall\ u \in U$.

\subsection{Example Applications}
We first demonstrate our results by extending the culminating finite example of \cite{RCRDSA} and giving a more complete picture. {\color{Red}As Hilbert spaces can be considered as approximation spaces, in our final example we consider the projective Hilbert space of complex numbers, $\mathbb{C}$. We show, among other things, the power set of the set of pure states and is a complete, atomistic Boolean algebra. We further suggest other Quantum relevant applications that may be useful, such as Hilbert spaces of operators. Finally, we touch upon references that are made to \textit{the statistical duality of states and effects} and some of the results stated in \cite{quant}. Those relationships between states and effects and probability measures also appear to have merit for further investigation.} We will continue to explore these concepts as we further develop this note.


\newpage
\begin{thebibliography}{9}
   \bibitem{Pawlak}
     Pawlak Z (1982) \emph{Rough sets.} Int J Comput Inf Sci 11:341–356
     
   \bibitem{RSTppf}
     Skowron, A., Dutta, S. \emph{Rough sets: past, present, and future.} Nat Comput 17, 855–876 (2018). https://doi.org/10.1007/s11047-018-9700-3
   
   \bibitem{CBgt}  
   Clouse D., Burke D. (2018) \emph{Approximating Power Indices to Assess Cybersecurity Criticality.} In: Bushnell L., Poovendran R., Başar T. (eds) Decision and Game Theory for Security. GameSec 2018. 
     
   \bibitem{RCRDSA}
     A. R. J. Srikanth, R. V. G. Ravi Kumar, G. V. S. R. Deekshitulu \emph{Core Rough Algebras and its Connection with Core Regular Double Stone Algebra} Journal of the Indian Math. Soc. Vol. 86, Nos. (3-4) (2019), 325-340
 
   \bibitem{BS}
     S. Burris, H. Sankappanavar \emph{A Course in Universal Algebra}, publicly released 2012, www.math.uwaterloo.ca/$\sim$ snburris

   \bibitem{CRDSA}
     A R J Srikanth, R V G Ravi Kumar \emph{Centre of Core Regular Double Stone Algbera}, European Journal of Pure and Applied Mathematics, Vol. 10, No. 4, 2017, pgs. 717-729

   \bibitem{pomykala}
     J. Pomykala and J. A. Pomykala, \emph{The Stone algebra of rough sets}, Bull. Polish Acad. Sci. Math.
36 (1988), 495–508.

   \bibitem{gamesec}
     Clouse D., Burke D. (2018) \emph{Approximating Power Indices to Assess Cybersecurity Criticality.} In: Bushnell L., Poovendran R., Başar T. (eds) Decision and Game Theory for Security. GameSec 2018. Lecture Notes in Computer Science, vol 11199. Springer, Cham. https://doi.org/10.1007/978-3-030-01554-1\_20
 
   \bibitem{nearlybool}Clouse D. \emph{The Nearly Boolean Nature of Core Regular Double Stone Algebras, CRDSA \tiny{(Ternary Set Partitions, CRDSA, Embeddings and Dual Equivalences)}} Cornell arXiv preprint https://arxiv.org/abs/1803.08313 

   \bibitem{comer}
     S. Comer \emph{On Connections Between Information Systems, Rough Sets and Algebraic Logic}, Algebraic Methods in Logic and Computer Science Banach Center Publicaions, Vol. 28, Institute of Mathematics Polish Academy of Sciences Warszawa, 1993

   \bibitem{duns}
     I. Duntsch, \emph{A logic for rough sets.}  Theoret. Comput. Sci., 179(1-2):427-436, 1997.

   \bibitem{dai}
     J.H. Dai, \emph{Rough 3-valued algebras}  Information Sciences, 178(8):1986-1996, 2008. http://www.sciencedirect.com/science/article/pii/S0020025507005361
     
   \bibitem{quant} 
     Cassinelli G, Lahti P. 2016 \emph{An axiomatic basis for quantum mechanics.} Found. Phys. 46, 1341–1373. (doi:10.1007/s10701-016-0022-y)
     
   \bibitem{quantc}
     Cassinelli G. and Lahti P. 2017 \emph{Quantum mechanics: why complex Hilbert space?} Phil. Trans. R. Soc. A.3752016039320160393 http://doi.org/10.1098/rsta.2016.0393 
\end{thebibliography}
\end{document}